\documentclass[12pt]{amsart}
\usepackage[margin = 1in]{geometry}
\usepackage{amsmath,amssymb,amsthm}
\usepackage[dvipsnames]{xcolor}
\usepackage{ulem}  
\usepackage{graphicx}
\usepackage{tikz,verbatim}

\usepackage{framed}

\newtheorem{theorem}{Theorem}

\theoremstyle{definition}
\newtheorem{definition}[theorem]{Definition}
\newtheorem{example}[theorem]{Example}

\title{Decomposable sparse polynomial systems}
\author[T.~Brysiewicz]{Taylor Brysiewicz} 
\address{T.~Brysiewicz\\ 
         Max Planck Institute for Mathematics in the Sciences\\ 
         Inselstr. 22, 04103\\ 
         Leipzig, Germany} 
\email{taylorbrysiewicz@gmail.com} 
\urladdr{https://sites.google.com/view/taylorbrysiewicz} 
\author[J.~I.~Rodriguez]{Jose Israel Rodriguez}
\address{J.~I.~Rodriguez\\Department of Mathematics\\
         University of Wisconsin\\
         Madison, WI 53706\\         
         USA}
\email{Jose@math.wisc.edu}
\urladdr{http://www.math.wisc.edu/~jose/}
\author[F.~Sottile]{Frank Sottile} 
\address{F.~Sottile\\ 
         Department of Mathematics\\ 
         Texas A\&M University\\ 
         College Station\\ 
         Texas \ 77843\\ 
         USA} 
\email{sottile@math.tamu.edu} 
\urladdr{http://www.math.tamu.edu/~sottile} 
\author[T.~Yahl]{Thomas Yahl} 
\address{T.~Yahl\\ 
         Department of Mathematics\\ 
         Texas A\&M University\\ 
         College Station\\ 
         Texas \ 77843\\ 
         USA} 
\email{thomasjyahl@math.tamu.edu} 
\urladdr{http://www.math.tamu.edu/~thomasjyahl} 
\thanks{Research of Sottile and Yahl supported by grant 636314 from the Simons Foundation.}

\newcommand{\CC}{\mathbb{C}}
\newcommand{\RR}{\mathbb{R}}
\newcommand{\ZZ}{\mathbb{Z}}
\newcommand{\calA}{\mathcal{A}}
\newcommand{\Adot}{{\calA_\bullet}}

\DeclareMathOperator{\MV}{{\rm MV}}

\newcommand{\defcolor}[1]{{\color{RoyalBlue}#1}}
\newcommand{\demph}[1]{\defcolor{{\sl #1}}}

\begin{document}

\begin{abstract}
 The Macaulay2 package  \texttt{DecomposableSparseSystems} implements methods for
 studying and numerically solving decomposable sparse polynomial systems.
 We describe the structure of decomposable sparse systems and explain how the methods in
 this package may be used to exploit this structure, with examples.  
\end{abstract}

\maketitle

\section{Introduction}

Am\'endola and Rodriguez~\cite{AmendolaRodriguez} gave numerical methods to efficiently solve systems of sparse polynomial
equations in a family, when that family is decomposable (Definition~\ref{D:decomposable}).
A consequence of Esterov's study of Galois groups of systems of sparse polynomial equations~\cite{Esterov} is that for
sparse systems, recognizing and computing a decomposition is algorithmic.
Solving a decomposable sparse system reduces to solving two smaller sparse polynomial systems.
In~\cite{SDSS}, we presented algorithms to detect and compute such decompositions, and a recursive algorithm 
exploiting decomposability for solving a decomposable sparse polynomial system using numerical homotopy continuation.

The Macaulay2 package \texttt{DecomposableSparseSystems} implements methods for 
decomposable sparse polynomial systems.
These include methods to detect decomposability, to compute a decomposition, and a recursive procedure to compute
numerical solutions to a given decomposable sparse system.
Detection and computation of decompositions uses integer linear algebra, including computing a Smith normal form and
the corresponding monomial changes of variables.
Numerical homotopy continuation is used to compute solutions.
When no further decompositions are possible, the algorithm solves multivariate systems using numerical software chosen by
the user (default: \texttt{PHCPACK}~\cite{V99}), and solves univariate polynomials using companion matrices.

Using the methods in \texttt{DecomposableSparseSystems} to solve a decomposable system allows for quicker solving and
more accurate solution counts than calling other solvers.
One reason is that after each decomposition, the child systems always involve either fewer variables, or polynomials of smaller degree. 
The cost of the methods  in \texttt{DecomposableSparseSystems} is low as they rely only on linear algebra and
numerical homotopy algorithms.

\section{Decomposable Sparse Polynomial Systems}

A \demph{branched cover} is a dominant map $\pi\colon X\to Y$ of irreducible varieties $X$ and $Y$ of the same dimension.
There is a number $d$ (the \demph{degree} of $\pi$) and an open dense subset $V$ of $Y$ such that $\pi^{-1}(v)$ consists of
$d$ points for $v\in V$. 
When $d>1$, the branched cover is \demph{nontrivial}.

\begin{definition}\label{D:decomposable}
 A branched cover $\pi\colon X\to Y$ is \demph{decomposable} if it is a composition of nontrivial branched covers.
 That is, if there is a dense open subset $U\subset Y$ and a variety $Z$ such that $\pi^{-1}(U)\to U$ factors as
\[
  \pi^{-1}(U)\ \longrightarrow\ Z\ \longrightarrow\  U\,,
\]
 with each map a nontrivial branched cover.
\hfill$\diamond$
\end{definition}
 In general it is not easy to determine if a branched cover is decomposable, or even to compute a decomposition for a
 decomposable branched cover.
 (See~\cite[Section~5.4]{AmendolaRodriguez} and~\cite[Section~1.2]{SDSS} for examples and a discussion.)

An integer vector $\alpha\in\ZZ^n$ is the exponent of a (Laurent) monomial
$\defcolor{x^\alpha} := x_1^{\alpha_1}\dotsb x_n^{\alpha_n}$.
A (complex) linear combination of monomials  $\sum c_\alpha x^\alpha$ is a (Laurent) polynomial.
Monomials are multiplicative maps $(\CC^\times)^n\to\CC^\times$ and polynomials are maps $(\CC^\times)^n\to\CC$.
For a finite set $\calA\subset\ZZ^n$ of exponents, the set of all polynomials whose monomials have exponents contained in $\calA$
(have \demph{support} $\calA$) forms the vector space \defcolor{$\CC^\calA$}. 
Given a list $\Adot = (\calA_1,\dotsc,\calA_n)$ of finite subsets of $\ZZ^n$, write \defcolor{$\CC^\Adot$} for the
vector space $\CC^{\calA_1}\oplus\dotsb\oplus\CC^{\calA_n}$ of lists $F=(f_1,\dotsc,f_n)$ of polynomials with $f_i$ having
support $\calA_i$.
Such a list $F\in\CC^\Adot$ is a function $F\colon(\CC^\times)^n\to\CC^n$,
and $F=0$ is a system of sparse polynomials with support $\Adot$ whose solutions are $F^{-1}(0)$.

\begin{example}\label{Ex:First}
Let $\Adot = (\calA_1,\calA_2)$ be the pair of supports in $\ZZ^2$ illustrated in Figure~\ref{Fig:one}.
\begin{figure}[htb]
\centering

\begin{picture}(85,71)(-14,0)
    \put(0,0){\includegraphics{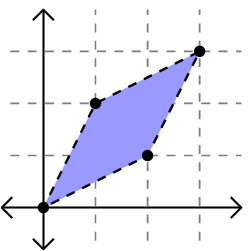}}
    \put(-14,19){$\calA_1$}
    \end{picture}
\qquad
\begin{picture}(100,71)(-14,0)
    \put(0,0){\includegraphics{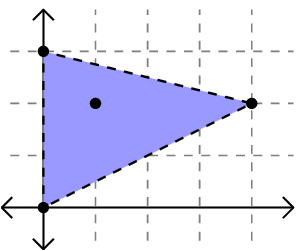}}
    \put(-14,19){$\calA_2$}
    \end{picture}

\caption{A pair of supports.}
\label{Fig:one}
\end{figure}
The corresponding vector spaces of polynomials are
\begin{eqnarray*}
   \CC^{\calA_1} &=&
           \left\{a_1+a_2xy^2+a_3x^2y+a_4x^3y^3\mid a_i\in\CC\right\}\,,\\
  \CC^{A_2}  &=& \left\{b_1+b_2y^3+b_3xy^2+b_4x^4y^2 \mid b_j\in\CC\right\}\,,
\end{eqnarray*}
and $\CC^\Adot$ is the space of systems of the form 
\[
   F\ =\  
\begin{pmatrix}
a_1+a_2xy^2+a_3x^2y+a_4x^3y^3 \\
b_1+b_2y^3+b_3xy^2+b_4x^4y^2
\end{pmatrix}\,,
 \quad
a_i,b_j\ \in\ \CC\,.
\]
In \texttt{DecomposableSparseSystems}, the family $\CC^\Adot$ is encoded by a list of matrices whose column vectors are
the exponent vectors of each polynomial. Given a system $F\in\CC^\Adot$,
these data can be
extracted from a given system via the Macaulay2 function \texttt{exponents}. \hfill$\diamond$ 
\end{example}

The Bernstein-Kushnirenko Theorem~\cite{Bernstein} provides a sharp {upper} bound on the number of solutions to a
system of sparse polynomials. 
Denote the convex hull of a set $\calA\subseteq\RR^n$ by \defcolor{$\text{conv}(\calA)$}.
Given a list of supports $\Adot=(\calA_1,\dotsc,\calA_n)$, let \defcolor{$\MV(\Adot)$} be the 
mixed volume of the list $(\text{conv}(\calA_1),\dotsc,\text{conv}(\calA_n))$.

\begin{theorem}[Bernstein-Kushnirenko]
  Let $\Adot$ be a list of $n$ finite subsets of $\ZZ^n$.
  For $F\in\CC^\Adot$, the number of isolated solutions in $(\CC^\times)^n$ to the system $F=0$ is bounded above by
  $\MV(\Adot)$ and this bound is achieved for $F$ lying in a dense, open subset of $\CC^\Adot$. 
\end{theorem}

Define $\defcolor{X_\Adot}\subset(\CC^\times)^n\times \CC^\Adot$ to be the set of pairs $(x,F)$ such that $F(x)=0$.
For $F\in\CC^\Adot$, the fiber $\pi^{-1}(F)$ of the map $\pi\colon X_\Adot\to\CC^\Adot$ consists of solutions to $F=0$.
By the Bernstein-Kushnirenko Theorem the map $\pi$ has degree $\MV(\Adot)$.
When $\MV(\Adot)\ge 1$, it is a branched cover.
When the branched cover $\pi:X_\Adot\to\CC^\Adot$ is decomposable, we say the sparse system $F\in\CC^\Adot$ is
decomposable.
Decomposability depends only on the support $\Adot$ of a system.

There are two transparent ways for a sparse system to decompose.

\subsubsection*{Lacunary}
A system $F\in\CC^\Adot$ is \demph{lacunary} if there is a surjective monomial map
$\Phi\colon(\CC^\times)^n\to(\CC^\times)^n$ such that $F = G\circ\Phi$ for some sparse polynomial system $G$.
We require that $\Phi$ be nontrivial in that its kernel is not the identity subgroup.
A lacunary system $F=G\circ\Phi=0$ can be solved by computing solutions, $z_1,\dotsc,z_d$, to the system $G=0$ and then
computing the fibres $\Phi^{-1}(z_1),\dotsc,\Phi^{-1}(z_d)$.
In appropriate coordinates, $\Phi$ is diagonal, and $\Phi^{-1}(z)$ is obtained by extracting roots of the components
of $z$.

\begin{example}
\label{Ex:lacunary}
Consider the following system with support from Example~\ref{Ex:First}.
 \[
   F(x,y) \ =\ 
   \begin{pmatrix}
   1-2xy^2+3x^2y-4x^3y^3\\
   2+3y^3+5xy^2+7x^4y^2
   \end{pmatrix}
   \ =\
   \begin{pmatrix} 0\\ 0 \end{pmatrix}
 \]
It is lacunary as it is the composition of the following maps.
\[
  G(s,t)\ =\ 
   \begin{pmatrix}
   1-2st^2+3st-4s^2t^3\\
   2+3st^3+5st^2+7s^2t^2
 \end{pmatrix}\,,
  \quad
  \Phi(x,y)\ =\ (x^3,x^{-1}y)\,.
\]
This can be detected via the methods in \texttt{DecomposableSparseSystems}.
\begin{leftbar}
\begin{verbatim}
i1 : R = CC[x,y];

i2 : F = {1-2*x*y^2+3*x^2*y-4*x^3*y^3,2+3*y^3+5*x*y^2+7*x^4*y^2};

i3 : isLacunary F
o3 = true
\end{verbatim}
\end{leftbar}
The method \texttt{isLacunary} extracts the set of supports of the system and computes the Smith normal form of a matrix
associated to these supports to determine whether the system is lacunary.\hfill$\diamond$
\end{example}

\subsubsection*{Triangular}
A system $F\in\CC^\Adot$ is \demph{triangular} if  
 there exists $k<n$ so that after a monomial change of variables, the system $F$ has the form  
\[
   F\ =\ (F_1(x_1,\dotsc,x_k),\dotsc,F_k(x_1,\dotsc,x_k),F_{k+1}(x_1,\dotsc,x_n),\dotsc,F_n(x_1,\dotsc,x_n))\,.
\]
Solutions to triangular systems are computed by first computing the solutions $z_1,\dotsc,z_d$ of the square subsystem
$(F_1,\dotsc,F_k)=0$.
A residual system is obtained by substituting $z_1$ into the original system for the first $k$ variables,
$F_2(z_1,x_{k+1},\dotsc,x_n)$.
Solutions to the original system are obtained by solving the residual system and then applying a homotopy
{algorithm} as described in~\cite{SDSS}. 

\begin{example}
Consider the system 
\[
  F(x,y,z)\ =\ 
 \begin{pmatrix}
   y^2-2x+3x^2y\\
   2+3x^2y+5x^4y^2
  \end{pmatrix}
  \ =\ 
\begin{pmatrix} 0\\ 0 \end{pmatrix}\,.
\]
Figure~\ref{Fig:two} shows the supports.
\begin{figure}[htb]
\centering

\includegraphics{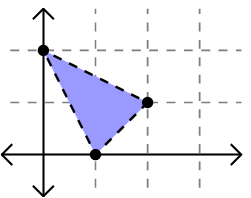}
\qquad
\includegraphics{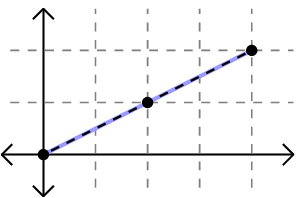}

\caption{Triangular support.}
\label{Fig:two}
\end{figure}
This system is triangular as the second polynomial is quadratic in the monomial $x^2y$.
The method \texttt{isTriangular} detects this subsystem.
\begin{leftbar}
\begin{verbatim}
i4 : F = {y^2-2*x+3*x^2*y,2+3*x^2*y+5*x^4*y^2};

i5 : isTriangular F
o5 = true
\end{verbatim}
\end{leftbar}\vspace{-15pt}
\hfill$\diamond$
\end{example}
A consequence of Esterov's study of Galois groups of sparse polynomial systems~\cite{Esterov} and Pirola and
Schlesinger's result that a branched cover is decomposable if and only if its Galois group is
imprimitive~\cite{PirolaSchlesinger} is that a sparse polynomial system is decomposable if and only if it is either
lacunary or triangular.
In each case, the solutions to original system are computed via solutions to simpler systems.
The methods in \texttt{DecomposableSparseSystems} iteratively decompose these sparse polynomial systems to efficiently
solve them.  

%
\section{Main method: solveDecomposableSystem}
The main method implemented in the package
\texttt{DecomposableSparseSystems} is named
\texttt{solveDecomposableSystem} and 
this implements Algorithm 9 in~\cite{SDSS}.
It takes as input a sparse polynomial system $F \in \mathbb{C}^{\Adot}$ and outputs all solutions
to $F=0$ in the algebraic torus.
 It recursively checks whether or not the input sparse polynomial system is decomposable, computes the decomposition, and then
calls itself on each portion of the decomposition. When the input is not decomposable it solves multivariate polynomial systems with the numerical solver
given by the option  \texttt{Software} (default: \texttt{PHCPACK}) and it solves univariate polynomial systems using companion matrices.
For complete details see \cite[Section~3.1]{SDSS}.

\subsection{Using the main method}
\begin{figure}[htb]
\centering

\begin{picture}(105,80)
    \put(0,0){\includegraphics{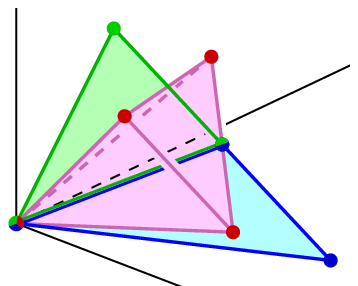}}
    \put(77,40){$\Adot$}
    \put(35,0){\scriptsize$x$}    \put(90,61){\scriptsize$y$}    \put(11,72){\scriptsize$z$}
    \end{picture}
\qquad
\begin{picture}(105,80)
    \put(0,0){\includegraphics{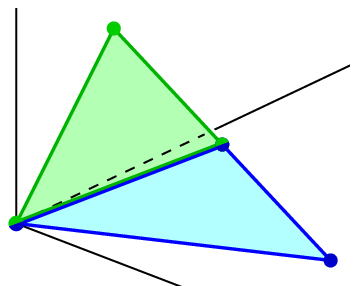}}
    \put(55,62){$\calA_3$}
    \put(90,25){$\calA_1$}
    \put(35,0){\scriptsize$x$}    \put(90,61){\scriptsize$y$}    \put(11,72){\scriptsize$z$}
    \end{picture}
\qquad
\begin{picture}(105,80)
    \put(0,0){\includegraphics{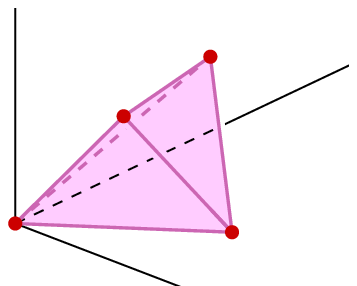}}
    \put(75,30){$\calA_2$}
    \put(35,0){\scriptsize$x$}    \put(90,61){\scriptsize$y$}    \put(11,72){\scriptsize$z$}
    \end{picture}

\caption{Support of $F$.}
\label{Fig:triple}
\end{figure}

Consider the system
\begin{align*}
F\ =\  
\begin{pmatrix}
2+xyz-x^2y\\
4-y^2z+2xz^2-3x^2z\\
1-yz^2-3xyz
\end{pmatrix}
\ =\  
\begin{pmatrix}
0\\
0\\
0
\end{pmatrix}.
\end{align*}
This system is supported on the triple $\Adot=(\calA_1,\calA_2,\calA_3)$ shown in Figure~\ref{Fig:triple}.

The method \texttt{isDecomposable} determines that this system is decomposable.
In particular, it is triangular with a subsystem indexed by the first and third polynomials.
This can be observed in the figure as the span of the supports $\calA_1$ and $\calA_3$ are coplanar.
It is also lacunary, as the exponent vectors lie in the sublattice of $\ZZ^3$ of index 3 generated by
the columns of $\left(\begin{smallmatrix}1&1&2\\ 1&0&0 \\ 1&2&1\end{smallmatrix}\right)$. The solutions to $F=0$ are found via the main method, \texttt{solveDecomposableSystem}.
\begin{leftbar}
\begin{verbatim}
i6 : R = CC[x,y,z];

i7 : F = {2+x*y*z-x^2*y,4-y^2*z+2*x*z^2-3*x^2*z,1-y*z^2-3*x*y*z};

  -- True if and only if the sparse system F is decomposable.
i8 : isDecomposable F
o8 : true

  -- A list of numerical solutions to F=0.
i9 : S = solveDecomposableSystem F; 

  -- Evaluates F at the first numerical solution.
i10 : F/(f-> sub(f, matrix {S_0}))
o10 = {1.77636e-15, 4.44089e-16+1.4623e-16*ii, 4.66294e-15}
\end{verbatim}
\end{leftbar}

Our main method also accepts a two-argument input $\texttt{(A,C)}$ where 
$\texttt{A}$ is a list of matrices whose columns support a system of (Laurent) polynomial equations, and $\texttt{C}$ is a list, whose $i$-th entry is the list of coefficients for the $i$-th polynomial equation.
We demonstrate some of the other types of inputs here, and leave details to the documentation.
\begin{leftbar}
\begin{verbatim}
i11 : (A,C) = (F/exponents/matrix/transpose,
                F/coefficients/last/entries/flatten);

i12 : S = solveDecomposableSystem (A,C);

  -- Expected timing for solving a specific system.
i13 : benchmark "solveDecomposableSystem(A,C)"; 
o13 = .0605920270512821

  -- Expected timing for solving a random system with support A.
i14 : benchmark "solveDecomposableSystem(A, )"; 
o14 = .0558867168108108
\end{verbatim}
\end{leftbar}

\subsection{Options for the main method}
Numerical in nature, \texttt{solveDecomposableSystem} features a variety of options for the user.
The option \texttt{Software} (default: \texttt{PHCPACK}) dictates which numerical solver is used to solve
multivariate sparse systems which are not decomposable.
The method \texttt{solveDecomposableSystem} removes solutions having any coordinate which is numerically zero up to
\texttt{Tolerance} (default: $10^{-5}$) throughout the computation.
Having this tolerance is necessary, as our methods are for Laurent polynomials with
solutions in the complex torus $(\CC^\times)^n$, while the solvers we call may return
solutions in $\CC^n$ that are not in the~torus.

When set to \texttt{1}, the option \texttt{Verify} (default: \texttt{0})
significantly increases the probability that \texttt{solveDecomposableSystem} computes the correct number of solutions. 

It does this by checking that \texttt{Software} computes $\MV(\Adot)$ solutions to any system $F$ with support $\Adot$, where $\MV(\Adot)$ is probabilistically  determined using \texttt{mixedVolume} in the package 
\texttt{Polyhedra}~\cite{Polyhedra}. 
If the mixed volume according  to \texttt{Polyhedra} and the number of solutions do not
agree, then the missing solutions are searched for using techniques related to those in 
\texttt{MonodromySolver} \cite{MonodromySolver}. 
Lastly, we allow the user to compute the solutions to $F$ by first solving an internally generated random instance and
then using that in a parameter homotopy~\cite{LSY1989} to solve $F$ 
by setting $\texttt{Strategy}$ to $\texttt{FromGeneric}$.
 We conclude by using the options \texttt{Verify} and  \texttt{Strategy} on an example with 6000 solutions.
\begin{leftbar}
\begin{verbatim}
i15 : A = <<< omitted, see example from Section 4 in [4] with
              i_1=(2,0,0,2,0)
              i_2=(4,4,2,2,2)
              j_1=(0,2,0,1,3)
              j_2=(0,0,1,0,2)
              >>;

  -- A has five supports, print the first one
i16 : print(length A,  A_0) 
(5, | 0 2 4 4 6 |)
    | 0 0 0 4 4 |
    | 0 0 0 2 2 |
    | 0 2 4 2 4 |
    | 0 0 0 2 2 |

i17 : elapsedTime (F,S) = solveDecomposableSystem(A,,Verify=>1);
     -- 8.93938 seconds elapsed

i18 : elapsedTime S' = solveDecomposableSystem(F,Strategy=>FromGeneric);
     -- 29.0802 seconds elapsed

i19 : print(#S,#S')
o19 = (6000, 6000)
\end{verbatim}
\end{leftbar}

\bibliographystyle{abbrv}
\bibliography{jsag}

\end{document}